\documentstyle[epic,eepic]{article}
\begin{document}
\title{ 
Strings, paths, and standard tableaux
}

\author{Srinandan Dasmahapatra\\
       \it Department of Mathematics,\\
       \it City University,\\
       \it London, EC1V 0HB, UK\\
       \\
       Omar Foda\\
       \it Instituut voor Theoretische Fysica,\\
       \it Universiteit Utrecht,\\
       \it Utrecht 3508 TA, The Netherlands \thanks{ 
       Permanent address: Department of Mathematics, 
       The University of Melbourne, Parkville, Victoria 
       3052, Australia}}
       
\maketitle

\begin{abstract}
For the vacuum sectors of regime-III ABF models, we observe that 
two sets of combinatorial objects: the strings which parametrize  
the row-to-row transfer matrix eigenvectors, and the paths which 
parametrize the corner transfer matrix eigenvectors, can both be 
expressed  in  terms  of  the  same  set  of  standard  tableaux. 
Furthermore, the momenta of the strings,  the  energies  of  the  
paths,  and the charges of the tableaux are such that there is a 
weight-preserving bijection between the two sets of eigenvectors, 
wherein the tableaux play an interpolating role.  This bijection 
is  so  natural, that  we  conjecture  that it exists in general.
\end{abstract} 

\newtheorem{theo}{Theorem}      
\newtheorem{cor}{Corollary}[theo]   
\newtheorem{de}{Definition}     
\newtheorem{pr}{Proposition} 
\newtheorem{co}{Corollary}[pr]  
\newtheorem{rem}{Remark} 
\newtheorem{ex}{Example}

\def\t(#1){\tilde{#1}}
\def\add{{\bf add}}
\def\sub{{\bf subtract}}
\def\cc{{\bf c}}
\def\pp{{\bf p}}
\def\ho{\hat{0}}
\def\hi{\hat{1}}

\def\nsection#1{\setcounter{equation}{0}\section{#1}}
\newtheorem{definition}{Definition}
\newtheorem{theorem}{Theorem}
\newtheorem{remark}{Remark}
\newtheorem{lemma}{Lemma}
\newtheorem{proposition}{Proposition}
\newtheorem{corollary}{Corollary}
\newtheorem{example}{Example}

\def\La{\Lambda}
\def\la{\lambda}
\def\Z{\:\mbox{\sf Z} \hspace{-0.82em} \mbox{\sf Z}\,}
\def\Zs{\mbox{\scriptsize \sf Z}  \! \! \mbox{\scriptsize \sf Z}}
\def\Znn{\Z_{\ge0}}
\def\C{\rule[0.5pt]{0.2mm}{6pt}{\hspace{-3.6pt}{\rm C}}}
\def\ie{{\em i.e., }}
\def\muh{\widehat{\mu}}
\def\P{{\cal P}_L}
\def\G{{\cal G}_L}
\def\Pset{\P(\La_i+\La_j,\La_k)}
\def\Gset{\G(2 \La_0,\La_k)}
\def\Mult#1#2#3{\left[#1 \atop #2 \right]_{#3}}
\def\MMult#1#2#3{\left[\left[#1 \atop #2 \right]\right]_{#3}}
\def\ib{\,\mbox{i}\,}
\def\is{\,\mbox{\scriptsize i}\,}
\def\e{\mbox{e}}
\def\es{\mbox{\scriptsize e}}
\def\case#1#2{{\textstyle{#1\over #2}}}
\def\eps{\epsilon}
\def\mod{\;(\bmod 2)}

\def\fow{{\tiny [FOW]}}
\def\kkr{{\tiny [KKR]}}

\section{Introduction} 

In statistical mechanics,  physical quantities are computed by 
averaging dynamical variables over an ensemble of all possible
state  configurations.  Each  configuration  is  weighted by a 
statistical  (Boltzmann)  factor.  The  transfer  matrix  is a 
convenient  tool  to generate the set of all properly-weighted 
configurations.  In  exactly-solvable  statistical  models  on 
two-dimensional lattices, the corner,  and row-to-row transfer 
matrices   (CTM   and   RTM,  respectively)  are  particularly 
useful.\footnote{For  an  introduction  to  exactly   solvable 
lattice models and transfer matrices,  see  \cite{BaxterBook}.} 

The  purpose  of this work is to point out a relation  between 
the   set of {\it paths} which parametrize the eigenvectors of 
the CTM \cite{BaxterBook,ABF,DJMO}, and the set  of  {\em strings,} 
which  play  the  same  role  for  the  RTM.
\cite{BRone}.  The  fact that such a relation must exist is 
not  new.  In  fact,  in  our  point  of view, it lurks in the 
background   of    all  recent   works  on  {\it boson-fermion} 
character  identities.\footnote{For  a  survey  of  the  early 
literature on this topic, see \cite{earlyreview}.  For  a more 
recent,  and relatively complete compendium of references, see 
\cite{BMO,Schilling}.}  What  is  new  in  this  work  is  the 
establishment of an explicit 
weight-preserving bijection between the relevant combinatorial 
objects: The two sets of eigenvectors,  and their  eigenvalues,
in  which,  to our delight, the  ubiquitous  standard tableaux
play a central role.  In the following, we briefly outline the 
contents and logic of the paper.   

We  will  restrict  our  attention  to  the anti-ferromagnetic 
(regime III) ABF models \cite{ABF},  because  they  are simple, 
well-understood and for our purposes, representative of  the 
general  situation.  These  models  form  an  infinite  series, 
labelled by  an integer $\ell = 1,~2, \cdots$.  The $\ell = 2$ 
model  is  the  Ising model (the $\ell=1$ model being trivial). 
The $\ell \geq 3$ models are {\it 
multi-critical} $Z_2$ models.\footnote{For  an  explanation of  
what these words  mean,  and  further  details  regarding  the 
physical   interpretation   of   these  models,  we  refer  to 
\cite{Huse}.} 
Working  in  a  certain  model,  one can specify {\it boundary 
conditions},  both  on  the outer boundary of the lattice, and 
the  value  of  a single variable at the center of the lattice. 
Each  choice of boundary conditions corresponds to a choice of 
{\it sector}  in  the  space of allowed lattice configurations. 
For  the  sake  of  simplicity,  in this work, we restrict our 
attention  to  the  {\it vacuum  sectors}  of  the  ABF models. 
\footnote{What  we  call  the vacuum sector corresponds to the 
choice $r=s=1$, of the parameters of \cite{ABF}.}

The   CTM   eigenvectors  can  be  parametrized  in  terms  of 
one-dimensional configurations  called  {\it paths}, that will 
be  defined  in  \S 2.1.  The corresponding eigenvalues, which 
turn   out   to   be   non-negative   integers,   are   called 
{\it energies},  and can be used as weights on the paths.  The 
weighted sum over all paths that belong to a certain sector in 
a  certain  model  turns out to be the  character of a highest 
weight module (HWM) of an infinite dimensional algebra. In the 
case of the ABF models, they are HWM's  of  Virasoro  algebras.  
This  sum can be evaluated by solving the recurrence relations  
satisfied  by  the  generating  function of the weighted paths.  
The  resulting expressions are in the form of alternating-sign 
q-series. Such  alternating-sign expressions are refered to in 
the   physics  literature  as  {\it bosonic}.\footnote{For   a  
discussion  of  the  origin  of  this terminology, we refer to 
\cite{earlyreview}.}

In \S 2, we relate the first set of combinatorial objects that 
we  are  interested  in:  the paths, to standard tableaux, and 
define  weights  on the tableaux, or tableaux statistics, that 
are  directly  related to the weights on the paths. The reason 
we  aimed  at  relating the paths to the tableaux is that the 
second  class  of combinatorial objects that we are interested 
in:   the   strings or rigged configurations\footnote{We shall 
be using the terms strings and rigged configurations
interchangably in this paper.},  which   label   
the RTM eigenvectors, and with which we wish to associate the paths, 
are also known to be related to standard tableaux \cite{KKR}. 

In  \S  3,  we  recall the definition of rigged configurations:
partitions  that are equipped, or {\it rigged}, with a set of
integers.  Both  the partitions, and the corresponding rigging 
integers satisfy certain model-dependent, and sector-dependent 
conditions.  Next,  we  outline a celebrated bijection between 
the  rigged  configurations  and  standard tableaux \cite{KKR}, 
and  recall the relation between the rigging integers, and the 
tableaux statistics mentioned above \cite{KR}.  

Summing  over  the  rigged configurations, of a certain sector, 
in  a  certain  model, where each configuration is weighted by 
a  sum  involving the  rigging  integers,  one obtains yet another 
expression for the character of the  corresponding  HWM.  Such 
expressions differ from those computed using CTM paths in that 
they   turn   out   to  be  constant-sign  expressions.  These
expressions are often referred to,  in  the physics literature, 
as   {\it  quasi-particle},  or  {\it  fermionic}  expressions 
\cite{KM,DKMM}.\footnote{For    further    details,   also 
regarding the derivation of the strings, their physical 
interpretation, and in  particular,  the physical significance  
of  the  rigging integers, which are related to the {\em Takahashi  
numbers}, or quasi-particle momenta in the physics literature, 
we refer to \cite{earlyreview}.}

Now  that  we  can  express  both  sets  of  eigenvectors, and 
corresponding eigenvalues,  in terms of tableaux and tableaux 
statistics respectively,  we identify the two sets of tableaux 
in \S 4.  This amounts to a bijective proof of the polynomial 
analogue of a boson-fermion character identity for each of the 
ABF  models  under  consideration.  The  reason that we obtain
polynomial,  rather  than  q-series identities, is that we can
consider paths, and rigged configurations, with an upper bound 
$L$  on  the  size of the system. In the limit 
$L \rightarrow \infty$, we
obtain  yet  another  proof of the generalized Roger-Ramanujan 
identities   conjectured   in  \cite{Melzer},  and  proven  in 
\cite{Berkovich,FodaWarnaar,Warnaar}. 

Though  we work in the limited context of the vacuum sector of 
regime III of the ABF models, the bijections that we establish 
are  so  simple  that  we  conjecture  that they extend to the 
entire   spectrum.  Furthermore,  we  believe  that  analogous 
bijections exist in more general models. These conjectures are 
supported  by  the  larger  number of conjectured, and in many 
cases   proven,  Rogers-Ramanujan-type  identities,  that  are 
related to more general ABF sectors, and other models. 

\section{CTM objects: paths and tableaux} 

\subsection{Paths} 

The  weighted sum over all configurations of a certain lattice
model,  with fixed dynamical variables on the boundaries,  and 
at  the  center  of  the  lattice,  defines  {\it a 
one-point function}. The one-point functions of the ABF models 
can be evaluated using Baxter's CTM method. 

Roughly speaking,  the CTM method reduces the physical problem 
of evaluating a one-point function of a lattice model,  to the 
mathematical  problem of evaluating the generating function of 
a set of {\it paths} on a segment of the weight lattice of the 
Lie algebra $sl(2)$. These paths are weighted in a certain way, 
and obey certain restrictions,  that  need not concern us here. 
All we need to know, is that the one-point functions depend on 
a  certain  {\it temperature-like}  parameter,  and  that  the
eigenvectors  of  the  CTM  consist  of linear combinations of 
paths  \cite{FodaMiwa}.  These  linear  combinations  can   be 
labeled  in  terms  of  a  single unique path \cite{Kashiwara}.  
In   the  zero  temperature  limit,  each  linear  combination  
reduces  to  the  single  path that labels it. The weight of a 
path  is  the  corresponding eigenvalue of the CTM in the zero 
temperature limit. 

In the case of vertex models based on affine algebras, such as 
the    six-vertex    model,   and   the   related   restricted 
solid-on-solid models, such as the ABF models, these one-point 
functions   can   be   expressed in terms of the characters of 
highest   weight   modules   of   affine and Virasoro algebras, 
respectively  \cite{DJMO}.  The  CTM  acts on the paths as the 
grading  operator, or derivation, of the corresponding algebra
\cite{FodaMiwa}. 

The  ABF  model characterized by the integer $\ell$ is defined 
in  terms  of {\it height variables} that live on the sites of 
a square lattice, and take values  in  the set of level $\ell$ 
dominant  integral  weights of the algebra $U_q\widehat{sl_2}$.  
A  path  is  a  sequence  of  heights where successive heights 
differ  by $\pm 1$.  Equivalently, for the fundamental weights 
${\Lambda_0,\Lambda_1}$, we introduce the set ${\cal P}_{\ell}$ 
which is the set of level $\ell$ weights:
$$
{\cal P}_{\ell} := \{\Lambda=k\Lambda_0 
                 + (\ell-k)\Lambda_1| 0\leq k\leq \ell\},
$$
and define 
${\cal A}:=\{\hat{0},\hat{1}\}$ 
where 
$\ho:=(\Lambda_1-\Lambda_0)$ 
and 
$\hi:=(\Lambda_0-\Lambda_1)$.  
Finally, introduce the automorphism $\sigma$ of the weight 
lattice which  acts  on  the  fundamental  weights  by the 
${\bf Z}_2$ transformation: 
$\Lambda_i \mapsto \sigma(\Lambda_i)=\Lambda_{1-i}$.

\begin{de}[path]
A $\Lambda$-path is a sequence $p=(\mu_0,\mu_1,\ldots,\mu_L)$, 
$\mu_j\in{\cal P}_{\ell}$ such that 
$\mu_{j+1}-\mu_j\in{\cal A}$, and for $j>L$, 
$\mu_j=\mu_L+\sigma^{(j-L)}(\Lambda)$.  We denote the set of all
such paths of length $L$ as ${\bf P}_\ell^L$.
\end{de}

We can encode a path in terms of a sequence of zeroes and
ones as follows:

\begin{de}[0 - 1 sequences]
For a path $p=(\mu_0,\cdots,\mu_L)\in {\bf P}^L_{\ell}$ we define 
a sequence $\eta(p)=(\eta_0,\cdots,\eta_L)$, 
where $\eta_j = \mu_{j+1}-\mu_j$.
\end{de}

In each sector, there is a unique path, $\bar{p}$ for which the 
energy  is lower than that of any other path in the same sector: 
{\em  the  ground  state path}. In regime-III of the ABF models, 
which is what we are interested in, $\eta(\bar{p})$ is given by 
$\eta_j \equiv \widehat{(j+1)}$ (mod $2$).

\par\noindent
 
\begin{ex}[ground state path]
$\bar{p} \in {\bf P}_{\ell}$ (note, we do not specify $L$ here),  
\begin{eqnarray*}
\bar{p}&=&(\ell\La_0,(\ell-1)\La_0+\La_1,
\ell\La_0,(\ell-1)\La_0+\La_1,
\ell\La_0,(\ell-1)\La_0+\La_1,
\ell\La_0,\cdots) \\
\eta(\bar{p})&=&(\ho, \hi, \ho, \hi, \ho, \hi,\cdots)
\end{eqnarray*}
\end{ex}
\par\noindent 
\begin{ex}
A path $p^{(1)}\in{\bf P}^L_\ell$.
\begin{eqnarray*}
p^{(1)}&=&(\ell\La_0,(\ell-1)\La_0+\La_1,
(\ell-2)\La_0+2\La_1,(\ell-1)\La_0+\La_1,\\
& & \quad \ell\La_0,(\ell-1)\La_0+\La_1,
\ell\La_0)
\\
\eta(p^{(1)})&=&(\ho, \ho, \hi, \hi, \ho, \hi)
\end{eqnarray*}
\end{ex}
We can represent the paths of the above examples as follows:

\smallskip

\setlength{\unitlength}{0.01in}
\begin{picture}(360,80)(0,-10)
\dottedline{5}(120,5)(120,65)(0,65)
	(0,5)(120,5)
\dottedline{5}(0,45)(120,45)
\dottedline{5}(0,25)(120,25)
\dottedline{5}(20,65)(20,5)
\dottedline{5}(40,65)(40,5)
\dottedline{5}(60,65)(60,5)
\dottedline{5}(80,65)(80,5)
\dottedline{5}(100,65)(100,5)
\dottedline{5}(360,5)(360,65)(240,65)(240,5)(360,5)
\dottedline{5}(240,45)(360,45)
\dottedline{5}(240,25)(360,25)
\dottedline{5}(260,65)(260,5)
\dottedline{5}(280,65)(280,5)
\dottedline{5}(300,65)(300,5)
\dottedline{5}(320,65)(320,5)
\dottedline{5}(340,65)(340,5)
\drawline(0,5)(20,25)(40,5)(60,25)(80,5)(100,25)(120,5)
\drawline(240,5)(280,45)(320,5)(340,25)(360,5)
\put(160,0){\makebox(0,0)[lb]{\raisebox{0pt}[0pt][0pt]{\shortstack[l]
{{$\, \, \scriptstyle \ell \La_0$}}}}}
\put(130,20){\makebox(0,0)[lb]{\raisebox{0pt}[0pt][0pt]{\shortstack[l]
{{$\, \, \scriptstyle (\ell-1) \La_0 + \La_1$}}}}}
\put(130,40){\makebox(0,0)[lb]{\raisebox{0pt}[0pt][0pt]{\shortstack[l]
{{$\, \, \scriptstyle (\ell-2) \La_0 + 2 \La_1$}}}}}
\end{picture}

\subsubsection{Weighted paths}

Let $p$ be a path, and $\bar{p}$ the corresponding ground-state path 
in ${\bf P}_\ell^L$,  with integer sequences $\eta(p)=(\eta_0,\ldots,
\eta_L)$, and  $\eta(\bar{p}) = (\bar{\eta}_0, \ldots,\bar{\eta}_L)$, 
respectively.  We  can  assign  to  each  path  a weight, called the
{\em energy function} as follows.  Introduce a functional 
$H:p\rightarrow {\bf Z}$:
\begin{equation}
H(p) = \sum_{j=1}^L j \theta(\eta_{j}-\eta_{j -1}),
\end{equation}
where $\theta$ is the step function 
\begin{equation}
\theta(z) = \left\{
\begin{array}{ll}
0 \qquad & (z  >   0) \\
1        & (z \leq 0),
\end{array}
\right.
\label{step}
\end{equation}
and define the energy $E(p)$ of the path $p$ as 
\begin{equation}
E(p)=H(p)-H(\bar{p}).
\end{equation}

Notice that, by construction, the ground state paths $\bar{p}$ of the 
vacuum sectors,  that we are interested in, have $H(\bar{p})=0$.  All 
other paths in the same set ${\bf P}^L_{\ell}$ have $H(p) > 0$. 
Comparing any path $p \in {\bf P}^L_{\ell}$ to $\bar{p}$, one can see 
that local departures from $\bar{p}$ \lq cost energy\rq, and that the 
further away from the origin the departure occurs,  the higher is the 
cost. 

\subsubsection{Polynomial analogues of branching functions}

We define polynomial analogues $B_L(q)$, of the branching functions 
or Virasoro characters
for the vacuum modules of the cosets 
$\widehat{sl(2)_{\ell-1}} \times \widehat{sl(2)_1} / 
\widehat{sl(2)_{\ell}}$
as the generating function of the weighted paths in ${\bf P}^L_{\ell}$
\cite{JMMO},
\begin{equation}
B_L(q) = \sum_{p\in {\bf P}^L_{\ell}}
q^{E(p)}.
\end{equation}
The paths that belong to the vacuum sectors, that we are interested in, 
begin and end at $\ell\Lambda_0$.

\subsection{Paths as tableaux}

A Young {\em diagram} $\lambda=(\lambda_1,\lambda_2,\lambda_3,\ldots)$   
is  a  collection  of cells arranged in rows of length $\lambda_1 \geq
\lambda_2 \geq\lambda_3\ldots$. A {\em standard} Young  tableau  is  a  
Young  diagram  whose  cells  are  occupied  by integers that strictly 
increase  along  rows  and  along  the columns. In this subsection, we  
present  a  bijection  between  the paths and standard tableaux, which 
allows us to use the latter to represent the former.  

To  each  path $p$, we associate a standard tableau $T_p$ as follows: 
scan the sequence of integers $\eta(p)$ that label a path $p$, and 
place  the  integer  $j$  in  the  $(i+1)^{\scriptstyle th}$ row, if
and only if $\eta_j=\hat{i}$.

We can adhere to {\rm normal} standard tableaux in this paper (no skew 
shapes) since we always have the first integer in  $\eta(p)=\ho$,  and 
the sequence of integers encoding such paths between dominant integral 
weights are {\em lattice permutations}.  

\begin{de}[lattice permutation] Consider an alphabet $\alpha$ in $n$ 
letters $\{\alpha_1,$ $\alpha_2, \cdots, \alpha_n\}$.  
The  letters are totally ordered in the sense 
that $\alpha_i < \alpha_j$  if $i < j$. A lattice permutation is a 
word in $\alpha$ such that,  on  scanning the letters that compose 
it from left to right, $\alpha_i$ occurs at least as many times as 
$\alpha_j$, for $i < j$.
\end{de}

Note that for a given $p$,  the contribution of step $j$,  to the 
path energy, is {\em zero} if the integer $(j+1)$ lies {\em below} 
$j$ in $T_p$.  With  this  property  in  mind,  we  introduce two 
integer-valued  functionals  or  {\em statistics} on the tableaux: 
$\cc(T)$ and $\pp(T)$.

\subsubsection{The charge of a tableaux, $\cc(T)$, and the Thomas 
functional $\pp(T)$} 

The functional $\pp(T)$, introduced by Thomas \cite{thomas}, is 
the sum of the 
positive integers  $i<n$ such that $i+1$ lies in a column to the 
{\em right} of $i$. (See \cite{macdonald}, p. 243.)
Observe that for any entry $j$ in a standard 
tableau $T$, $j+1$ lies {\em either} to the right of $j$ {\em or} 
below   it.   Therefore,   $\pp(T_p)$   picks  up  the  positive 
contributions to $\theta(\eta_j-\eta_{j-1})$ for all $j$. 

The  charge  of  a  tableau, $\cc(T)$ was defined by Lascoux and 
Sch\"utzenberger \cite{L-Sch} for a tableau $T$ of shape $\lambda$ 
and weight $\mu$. (The weight of a tableau is a sequence 
$\{\mu_1,\mu_2,\ldots,\mu_n\}$, 
where  $\mu_i$  is the number of times $i$ occurs in the tableau.  
In  this  paper we are interested in standard tableau, {\em i.e.} 
the  case  $\mu=\{1, 1, \cdots, 1\}$.)  We shall follow 
\cite{macdonald}, 
page 242,  where $\cc(T)$ is defined for any tableau whose weight 
$\mu$ is a partition.

{}From a tableau $T$, read off the entries from right to left and 
from top to bottom to get a word $w(T)=a_1 a_2\ldots a_n$. Assign 
an {\em index} to each $a_i$ as follows.  Let  $1$ have index $0$. 
If the number $r$ has index $i$,  then  $r+1$ will have index $i$ 
if  it  is  to  the {\em right} of $r$, and $i+1$ if it is to the 
{\em left}.  The  charge  $c(w)$  of a word $w$ is the sum of its 
indices.  The  charge  of  a  tableau $T$ is the charge $c(w(T))$.

For example, for the tableau
\[
T=
\begin{array}{cccc}
1 & 2 & 4 & 8 \\
3 & 6 & 9 &   \\
5 & 7 &   &  
\end{array},
\]
the  word $w(T)=8 4 2 1 9 6 3 7 5$.  The indices of $1,2,\ldots,9$ 
are  $0,1,1,2,2,3,3,4,4$, respectively  and  $c(w(T)) = 20$.  Also, 
$\pp(T)=1+3+5+7=16$.

\subsubsection{The Sch\"utzenberger involution}

The  involution  we shall describe is a special case of a family of 
operations which are indexed by points $p\in Z\times Z$ occupied by 
an ordered set of objects on which they act injectively. These were 
defined by Sch\"utzenberger and are called {\em jeu de taquin}. For 
the  involution  $S$  on  standard  tableaux,  also  called {\em an 
evacuation},  the  point in question is the top left corner $(1,1)$, 
{\em e.g.},  the  cell  occupied  by  the  number  $1$ in the above 
example.  

The  operation $S$ is defined by the following set of moves.  First, 
remove  the  entry  in  cell $(1,1)$.  Move  the smaller of the two 
integers immediately below or to the right of it to occupy the cell 
$(1,1)$.  This  leaves a vacancy in the cell from which the integer 
was moved.  Then look for the smallest integer immediately below or 
to the right of this freshly vacated cell,  and move that one in to 
squat in it. This leaves a new cell empty. In general, whenever the 
cell $(i,j)$ is empty, and 
$t(i_1,j_1)= \mbox{min}(t(i+1,j),t(i,j+1))$, 
slide  the  integer  from  $(i_1,j_1)$ to $(i,j)$. And so on, until 
there  are  no  more  cells to the right and below the cell vacated.  
In  this last cell insert the number $n$ with parentheses around it, 
to  distinguish  it  from  the  $n$  that is already present in the 
tableau.  Now,  remove the (new) integer from the cell $(1,1)$, and 
repeat  this procedure, inserting $(n-1)$ in the last cell.  And so 
on,  until all the entries of the original tableau are deleted, and 
we have parentheses around the integers in the new standard tableau, 
which we then remove.

\begin{ex}[An evacuation]
In  the  tableau  $T$  above, remove the entry from the cell $(1,1)$ 
and start the evacuation:
\[
\begin{array}{cccc}
\cdot & 2 & 4 & 8 \\
3 & 6 & 9 &   \\
5 & 7 &   &  
\end{array}
\rightarrow
\begin{array}{cccc}
2 & \cdot & 4 & 8 \\
3 & 6 & 9 &   \\
5 & 7 &   &  
\end{array}
\rightarrow
\begin{array}{cccc}
2 & 4 & \cdot & 8 \\
3 & 6 & 9 &   \\
5 & 7 &   &  
\end{array}
\rightarrow
\begin{array}{cccc}
2 & 4 & 8 & \bullet\\
3 & 6 & 9 &   \\
5 & 7 &   &  
\end{array}.
\]
Replace the $\bullet$ by $(9)$. Iterating this one more time, we get
\[
\begin{array}{cccc}
\cdot & 4 & 8 & (9) \\
3 & 6 & 9 &   \\
5 & 7 &   &  
\end{array}
\rightarrow
\begin{array}{cccc}
3 & 4 & 8 & (9) \\
\cdot & 6 & 9 &   \\
5 & 7 &   &   
\end{array}
\rightarrow
\begin{array}{cccc}
3 & 4 & 8 & (9) \\
5 & 6 & 9 &   \\
\cdot & 7 &   &  
\end{array}
\rightarrow
\begin{array}{cccc}
3 & 4 & 8 & (9) \\
5 & 6 & 9 &   \\
7 & \bullet &   &  
\end{array}.
\]
After a couple of more steps, we have
\[
\begin{array}{cccc}
3 & 4 & 8 & (9) \\
5 & 6 & 9 &   \\
7 & (8) &   &  
\end{array}
\rightarrow
\begin{array}{cccc}
4 &  6 &  8 & (9)\\
5 &  9 & (7) &   \\
7 & (8) &   &  
\end{array}
\rightarrow
\begin{array}{cccc}
5  & 6 & 8 & (9)\\
7  & 9 & (7) &   \\
(6) & (8) &   &  
\end{array},
\]
and so on until 
\[
T^S=
\begin{array}{cccc}
1 & 3 & 5 & 9 \\
2 & 4 & 7 &   \\
6 & 8 &   &  
\end{array}.
\]
\end{ex}
\noindent It is possible to show that the following two properties 
in the above example:

\noindent {\it i)} $S$ is an involution, {\em i.e.} $S^2=1$, or
$(T^S)^S=T$.

\noindent {\it ii)} $\cc(T)=\pp(T^S)$.

\noindent Both of these properties are, in fact, true in general.
See \cite{Butler} and references therein.

\subsubsection{Energy of a path as the statistic {\bf p}}

The following is a key observation: Let a path 
$p=(\eta_1,\ldots,\eta_L)$ 
be encoded as a Young tableau $T_p$. The contributions to $H(p)$ come 
from  the  integers  $j$  for  which $(j+1)$ is to the right of it in 
$(T_p)$.  This, together with the definition of ${\bf p}(T_p)$, leads
us to

\begin{pr}
For a path $p=(\eta_1,\ldots,\eta_L)$ encoded as a Young tableau
$T_p$ with a total number of nodes $|T_p|=L$, 
$H(p)={\bf p}(T_p).$  The energy of the path $E(p)$ is then
$$E(p)={\bf p}(T_p)-{L\over 2}({L\over 2}-1).$$

\end{pr}

\noindent {\em Proof:} 
Observe that for a path $\bar{p}:=(0,1,0,1,\ldots,0,1)$, of length $L$, 
if the corresponding tableau is $\bar{T}$, then 
$${\bf p}(\bar{T})
=
\sum_{i=1}^{L/2-1}(2i)={L\over 2}({L\over 2}-1).$$  The proof follows
directly from the definition of $E(p)$.
{\hfill $\Box$}

\section{RTM objects: strings and tableaux}

In \cite{BRone} the eigenvalues of the row-to-row transfer matrix 
of  the  ABF  models  were  obtained in terms of solutions to the 
Bethe equations.  In  order  to  solve  the  Bethe equations, one
typically  starts  by  {\em assuming}  that  the  solutions  form 
{\em strings}  \cite{Bethe}:  they  form  clusters in the complex 
plane, where each cluster of roots has a common 
real part, and equally spaced imaginary  parts.  A  cluster  that
consists  of  $j$  of  elements  is called a string of length $j$.

It is standard \cite{Tak} to multiply the Bethe equations for the 
components of each string, and end up with equations for the real 
parts of the solutions.  The  logarithm of the multiplied form of 
the  Bethe  equations is then taken such that the (half-) integer 
branches  are {\em  distinct},  hence  referred  to as {\em 
fermionic}. 

Every  eigenvector of the RTM can be described in terms of such a 
set of (half-) integers,  which form convenient labels for describing 
the 
physics of these models. The characterization of roots in terms of the
length of the strings and the (half-)integer  branches  of  logarithms  
(also called Takahashi  numbers) is  called a
{\em rigged configuration} in \cite{KKR}. Let us express 
the   indexing   and  counting  of  states  in  the  language  of 
\cite{KKR},\cite{KR}.

For a given two rowed partition (or Young diagram)
$\lambda=(\lambda_1,\lambda_2), \lambda_1+\lambda_2=L$, 
we shall define another, $\nu, \nu\vdash \lambda_2$ so that
the rows of a partition $\nu$ represent strings, and the 
lengths of these rows are viewed as the lengths of these 
strings.  We shall require that the maximum length of 
strings allowed is $\ell$.
Let $m_j$ be the number of rows 
length $j$, and $Q_j$ be the number of cells in the first 
$j$ columns of diagram $\nu$.

The rows are labelled by integers 
$J_{j,\mu},\mu=1,\ldots,m_j,j=1,\ldots,\ell$,
$0\leq J_{j,\mu}\leq P_j$, and we require $P_j\geq 0$
for a rigged configuration to be {\em admissible}.
Unlike the fermionic Takahashi numbers, these integers 
$J_{j,\mu}$, however, are
allowed to repeat, {\em i.e.} they may be called 
{\em bosonic}.\footnote{For a
fixed number of particles, the labels bosonic or fermionic do not, 
however, mean very much.  It is only when the particle number changes
that the new occupation numbers for bosons or fermions or any other '-ons'
exhibit their statistical character.} 

The maximum allowed {\em bosonic} integer $P_j$
for any of the labels of the 
$m_j$ strings of length $j$ is given by $P_j=L-2Q_j$.
Note, $Q_j=\sum_{i=1}^{\ell}$min$(i,j)m_i$ and therefore,
\begin{eqnarray}
P_j & = & L-2\sum_{i=1}^j \sum_{k=i}^\ell m_k \\
    & = & (\lambda_1-\lambda_2)+2\sum_{k>j}^\ell (k-j)m_k.
\end{eqnarray}

\begin{de}[rigged configurations]
A configuration of partitions $\nu$ whose rows are
indexed by non-negative integers 
$J_{j,\mu}, \mu=1,\ldots,m_j$ 
no larger than $P_j$ for $1,\leq j\leq \ell$
where $P_j$ is calculated as above
is called a {\em rigged configuration}.  We shall
denote a rigged configuration by $(\nu,L,\{J\})$
\end{de}

The total number of such rigged configurations is given by 
\begin{equation}
\sum_{m_1,\ldots,m_\ell}\prod_{j+1}^\ell\biggl(
\begin{array}{c}
P_j+m_j\\
m_j
\end{array}\biggr)
\end{equation}

\subsubsection*{Remark}  The counting procedure above is identical to 
the one involving Takahashi integers $I_{j,\mu}$ (see\cite{das}), 
which are related to the (bosonic) integers $J_{j,\mu}$ by
$$I_{j,\mu}=J_{j,\mu}+\mu-
{1\over2}(P_j+m_j+1).
$$

\subsection{The quasi-particle momentum of a rigged configuration} 

While the sum of the integers gives the momentum of the eigenstate, 
fermionic character sums are generating functions for the momenta
of sectors of the excited states in the theory, chosen in an 
appropriate way.  These excited states are labelled by their 
quasi-particle content, and for the correct counting of quasi-particle 
momenta, it is necessary to assign the 
appropriate zero-momentum point.  This is done by subtracting off the 
{\em Fermi momentum}, and for our models, and the excitations counted
by the variable $P_k$, the edge is given by the $L$ dependent piece in 
equation .  We therefore have to take the following sum in order to
get the momentum of the excited states counted by $P_k$:
\[
P(\nu,\{I\})=\sum_{j=1}^{\ell-1} \sum_{\alpha=1}^{P_j} 
(I_{j,\alpha}-{L\over 4} \delta_{1j}).
\]
A useful way of viewing this sum for the momentum of the quasiparticle
states is to write it in terms of the variables $J_{j,\alpha}$
decompose it into a piece which has all the $J_{j,\alpha}=0$
and a sum over the integers $J_{j,\alpha}$.  That is,
\begin{equation}
\label{qp}
P(\nu,\{J\})={1\over 4}\sum_{i,j=1}^{\ell-1} P_i C_{i,j}P_j
+\sum_{\alpha=1}^{P_j}
J_{j,\alpha},
\end{equation}
where $C_{i,j}$ are the elements of the Cartan matrix for the $A_{l-1}$
root system.

\subsection{Rigged configurations and lattice permutations}

Before launching into the description of the bijective correspondence
between rigged configurations and tableaux, or equivalently,
lattice permutations, let us first prove the following proposition.

\begin{pr}
A rigged configuration is admissible  
if and only if the tableau associated to it is a standard tableau.
\end{pr}
\noindent {\em Proof:}
Let us assume that the tableau $T$ associated to a sequence of $0$'s and 
$1$'s is not standard, {\em i.e.,} $\lambda_1-\lambda_2<0$ for the shape
$\lambda=(\lambda_1,\lambda_2)$ of $T$.  This implies $P_\ell<0$.
Conversely, note that
$$P_{j_1}-P_{j_2}=2\sum_{k=j_1+1}^{j_2} (k-j_1)m_k + 
2\sum_{k=j_2+1}^\ell (j_2-j_1) m_k \geq 0, \,\, j_1<j_2.$$
This implies that if the configuration is not admissible, 
$P_j<0$ for some $j$ $\Longrightarrow P_\ell<0$ and therefore 
$\lambda_1 - \lambda_2<0$.
\hfill{$\Box$}

\subsubsection{From a lattice permutation to a rigged configuration}

A key idea, used repeatedly is that of a {\em special string} or
row of a rigged configuration.
\begin{de}[special string] 
A row in $\nu$ of length $i$ and label $\alpha, \alpha=1,\ldots,
m_i$ in a rigged configuration 
($\nu,L,\{J\}$) is called `special' if its rigging integer 
is maximal, {\em i.e.} $J_{i,\alpha}=P_i$.
\end{de}

{}From the sequence of integers $w=b_1 b_2 \ldots b_L$, we shall 
construct the rigged configuration that corresponds to it by an 
inductive process due to Kerov, Kirillov and Reshetikhin \cite{KKR}.

Let us assume that we have already constructed $(\nu,k-1,\{J\})$
from $w_1=b_1 b_2 \ldots b_{k-1}$.  Order the rows of equal length
so that the integers within each such block are weakly decreasing
from the top.  So, of the rows are numbered $1,\ldots,m_j$ from the
top of a block of $j$-strings, 
\begin{equation}
\label{dec}
P_j\geq J_{j,1}\geq \cdots \geq J_{j,m_j}\geq 0 \quad \forall j.
\end{equation}
We shall now describe how to get to
$(\t(\nu),k,\{\t(J)\})$, given the value of $b_k$.  There are two
possibilities: $b_k=0,1$.  For $b_k=0$, set $\t(\nu)=\nu$ and 
every rigging integer $\t(J)=J$ (we suppress the labels here).  The only
change is therefore, $k-1\mapsto k$, and therefore $P_j\mapsto P_j +1,
\forall j$.

If $b_k=1$, set $i^*=max\{i|J_{i,1}=P_i\}$ and $r^*=1+\sum_{i>i^*}m_i$.
Then 
$$\nu=(\nu_1,\ldots,\nu_{{\scriptstyle Q_1}})\mapsto\t(\nu)=
(\nu_1,\ldots,\nu_{r^*}+1,\ldots,\nu_{\scriptstyle Q_1}).$$
In words, we add a box to the highest special row. 

Next, we need to describe $\{\t(J)\}$.  We set $\t(J)=J$ for all 
rows of $\t(\nu)$ except the $r^{* \scriptstyle th}$ row, which 
we set equal to the maximum allowed.  Note that the following
changes occur to the maximum integers $P_j$:

\begin{equation}
P_j\mapsto \t(P)_j= \left\{
\begin{array}{ll}
P_j + 1,& (j\leq i^*)\\
P_j - 1,& (j>i^*).
\end{array} 
\right.
\label{P-add}
\end{equation}
Thus the new (changed) integer is set equal to 
$\t(P)_{i^*+1}=P_{i^*+1}-1$.  Once again, we re-order the rows of
length $i^*+1$ so that the integers satisfy (\ref{dec}).
{}From (\ref{P-add}) it is also clear that $\t(P)_j\geq 0\,\forall j$,
so that the new configuration is also admissible.

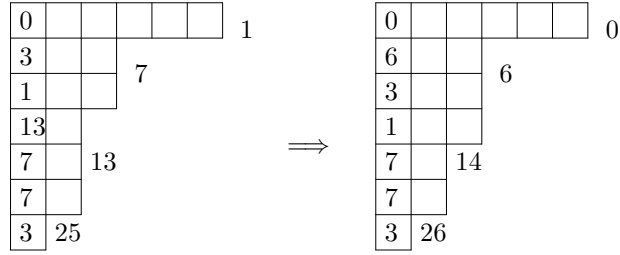
\begin{ex}[$b_k=1$]

\centerline{
\setlength{\unitlength}{0.00925in}
\begin{picture}(136,155)(0,-10)
\drawline(0,140)(120,140)
\drawline(0,80)(60,80)
\drawline(0,0)(20,0)
\drawline(0,20)(40,20)
\drawline(120,140)(120,120)(60,120)(60,80)
\drawline(40,140)(40,20)
\drawline(20,140)(20,20)
\drawline(0,140)(0,20)
\drawline(0,20)(0,0)
\drawline(20,20)(20,0)
\drawline(0,40)(40,40)
\drawline(0,60)(40,60)
\drawline(0,100)(40,100)
\drawline(0,120)(40,120)
\drawline(40,100)(60,100)
\drawline(40,120)(60,120)
\drawline(60,140)(60,120)
\drawline(80,140)(80,120)
\drawline(100,140)(100,120)
\put(130,120){\makebox(0,0)[lb]{\raisebox{0pt}[0pt][0pt]{\shortstack[l]{{\rm 1}}}}}
\put(70,95){\makebox(0,0)[lb]{\raisebox{0pt}[0pt][0pt]{\shortstack[l]{{\rm 7}}}}}
\put(45,45){\makebox(0,0)[lb]{\raisebox{0pt}[0pt][0pt]{\shortstack[l]{{\rm 13}}}}}
\put(25,5){\makebox(0,0)[lb]{\raisebox{0pt}[0pt][0pt]{\shortstack[l]{{\rm 25}}}}}
\put(5,125){\makebox(0,0)[lb]{\raisebox{0pt}[0pt][0pt]{\shortstack[l]{{\rm 0}}}}}
\put(5,105){\makebox(0,0)[lb]{\raisebox{0pt}[0pt][0pt]{\shortstack[l]{{\rm 3}}}}}
\put(5,85){\makebox(0,0)[lb]{\raisebox{0pt}[0pt][0pt]{\shortstack[l]{{\rm 1}}}}}
\put(5,65){\makebox(0,0)[lb]{\raisebox{0pt}[0pt][0pt]{\shortstack[l]{{\rm 13}}}}}
\put(5,45){\makebox(0,0)[lb]{\raisebox{0pt}[0pt][0pt]{\shortstack[l]{{\rm 7}}}}}
\put(5,25){\makebox(0,0)[lb]{\raisebox{0pt}[0pt][0pt]{\shortstack[l]{{\rm 7}}}}}
\put(5,5){\makebox(0,0)[lb]{\raisebox{0pt}[0pt][0pt]{\shortstack[l]{{\rm 3}}}}}
\end{picture}
\raisebox{.6in}{$\quad\Longrightarrow\quad$}
\setlength{\unitlength}{0.00925in}
\begin{picture}(136,155)(0,-10)
\drawline(0,140)(120,140)
\drawline(0,80)(60,80)
\drawline(0,0)(20,0)
\drawline(0,20)(40,20)
\drawline(120,140)(120,120)(60,120)(60,80)
\drawline(40,140)(40,20)
\drawline(20,140)(20,20)
\drawline(0,140)(0,20)
\drawline(0,20)(0,0)
\drawline(20,20)(20,0)
\drawline(0,40)(40,40)
\drawline(0,60)(40,60)
\drawline(0,100)(40,100)
\drawline(0,120)(40,120)
\drawline(40,100)(60,100)
\drawline(40,120)(60,120)
\drawline(60,140)(60,120)
\drawline(80,140)(80,120)
\drawline(100,140)(100,120)
\drawline(40,60)(60,60)(60,80)
\put(5,125){\makebox(0,0)[lb]{\raisebox{0pt}[0pt][0pt]{\shortstack[l]{{\rm 0}}}}}
\put(5,45){\makebox(0,0)[lb]{\raisebox{0pt}[0pt][0pt]{\shortstack[l]{{\rm 7}}}}}
\put(5,25){\makebox(0,0)[lb]{\raisebox{0pt}[0pt][0pt]{\shortstack[l]{{\rm 7}}}}}
\put(5,5){\makebox(0,0)[lb]{\raisebox{0pt}[0pt][0pt]{\shortstack[l]{{\rm 3}}}}}
\put(5,65){\makebox(0,0)[lb]{\raisebox{0pt}[0pt][0pt]{\shortstack[l]{{\rm 1}}}}}
\put(5,105){\makebox(0,0)[lb]{\raisebox{0pt}[0pt][0pt]{\shortstack[l]{{\rm 6}}}}}
\put(5,85){\makebox(0,0)[lb]{\raisebox{0pt}[0pt][0pt]{\shortstack[l]{{\rm 3}}}}}
\put(130,120){\makebox(0,0)[lb]{\raisebox{0pt}[0pt][0pt]{\shortstack[l]{{\rm 0}}}}}
\put(70,95){\makebox(0,0)[lb]{\raisebox{0pt}[0pt][0pt]{\shortstack[l]{{\rm 6}}}}}
\put(45,45){\makebox(0,0)[lb]{\raisebox{0pt}[0pt][0pt]{\shortstack[l]{{\rm 14}}}}}
\put(25,5){\makebox(0,0)[lb]{\raisebox{0pt}[0pt][0pt]{\shortstack[l]{{\rm 26}}}}}
\end{picture}
}

$(\nu,L=39)\Longrightarrow(\nu',40)$, for $b_{40}=1$.  A node is added
to the highest special row, and the resulting configuration has its rows
rearrangeed so that the integer riggings are weakly decreasing.
\end{ex}

\subsubsection{From a rigged configuration to a lattice permutation}

Here we describe a way of reading off $w_k$ from a 
rigged configuration ($\nu,k,\{J\}$), where $w_k$ is 
either a $0$ or a $1$, and modifying it to get
another configuration ($\bar{\nu},k-1,\{\bar{J}\}$),
a process called {\em ramification}.  We can perform
the same procedure on ($\bar{\nu},k-1,\{\bar{J}\}$), 
to extract $w_{k-1}$, so this is 
a recursive method for generating a sequence of $0$'s and $1$'s 
from a given rigged configuration ($\nu,L,\{J\}$).  

Here we describe the {\em ramification rules}.
First, we arrange the rows of equal length $j$ so that the integers are
in a weakly increasing order from the top:
\begin{equation}
\label{inc}
0\geq J_{j,1}\leq \cdots \leq J_{j,m_j}\leq P_j \quad \forall j.
\end{equation}
Set $i_*=min\{i|J_{i,m_i}=P_i\}$ and $r_*=\sum_{i\geq i_*}m_i$, and
also assign
$$\nu=(\nu_1,\ldots,\nu_{\scriptstyle Q_1})\mapsto\bar{\nu}=
(\nu_1,\ldots,\nu_{r_*}-1,\ldots,\nu_{\scriptstyle Q_1}).$$
In words, we delete a box from the lowest special row. 
		     
Next, we need to describe $\{\bar{J}\}$.  We set $\bar{J}=J$ for all 
rows of $\bar{\nu}$ except the $r_{*}^{\scriptstyle th}$ row, which 
we set equal to the maximum allowed.  Note that the following
changes occur to the maximum integers $P_j$:
\begin{equation}
\label{P-del}
P_j\mapsto \bar{P}_j=\left \{
\begin{array}{ll}
P_j - 1,& (j < i^*)\\
P_j + 1,& (j \geq i^*).
\end{array}
\right.
\end{equation}
Thus the new (changed) integer is set equal to 
$\t(P)_{i^*-1}=P_{i^*-1}-1$.  Once again, we re-order the rows of
length $i^*+1$ so that the integers satisfy (\ref{inc}).
{}From (\ref{P-del}) it is also clear that $\bar{P}_j\geq 0\,\forall j$,
so that the new configuration is also admissible.

We can now define $w_k:=0$ if $i_*=0$, {\em i.e.}, if there are no
special strings in $\nu$, or $w_k:=1$ if there is at least one such
special string ($i_*>0$).  

\section{The bijection between paths and rigged configurations}

The  starting point of our proposed bijection between the CTM paths, 
and  the RTM rigged configurations,  are the following observations: 
Once  we  know that the Bethe Ansatz solution of a model can be cast 
in terms of interacting quasi-particles \cite{earlyreview},  one can 
attempt  re-phrase  the  CTM solution of the same model in that same 
language.  The  ABF  model labelled by $\ell$ 
contains strings of length $1,\ldots,\ell$.  Looking  at  the  CTM 
paths of such a model, one can think of the {\em peaks} of different 
heights  as  representing the various types of strings. (One 
can think of the peaks of height 1 as {\em non-physical}:  they  are 
there  to  fill  the  {\em Dirac  ({\em or} Fermi sea}, but do not 
contribute to the 
total momentum of the physical states.)  

In this picture, each path encodes the {\em momentum} of a configuration 
of strings.  It is now 
clear that one can classify the set of all paths into sectors,  such 
that,  the  set of all paths in each sector share the same number of 
string-lengths, and therefore quasi-particles of each type,  
but  the different  paths in each sector correspond  to 
different excited states of the same set of quasi-particles. 

It turns out,  and it is easy to see by inspection, that each sector 
will  conatin  a  certain  unique  path,  with lowest possible total 
momentum.  We  wish  to  refer  to  this  as  {\em  a  minimal  path 
configuration}.  Furthermore,  since  we are working on finite paths, 
each  sector  will also contain a unique path, with largest possible 
total  momentum.  We  wish  to  refer  to  that  as {\em  a  maximal 
configuration}.  As we will see below, identical remarks can be made 
about the rigged configurations, and for a rigged configuration of a 
certain  set  of  strings,  one  can  define a minimal and a maximal 
configuration.  In order to identify the paths and the corresponding 
rigged configurations, we need to make sure that 
we do have a weight-preserving bijection.

We  proceed  in  two  steps: firstly, we show that, sector by sector, 
the  minimal  and  maximal  path  and  rigged  configurations are in
bijection. Secondly, we show that, sector by sector, the rest of the 
configurations are also in bijection with the other paths. 

\subsection{The minimal configurations}

Let $(\nu,L,\{J=0\})$ be a rigged configuration with a string content
$(m_1,m_2,\ldots,$ $m_l)$,  with all the rigging integers set to zero. 
Since some of the $m_i$ could be zero, we shall label the rows from 
$i=1$ (the first row) to $i=Q_1=m_l+m_{l-1}+\cdots+m_1$ (the bottom 
row) and let $\nu_i$ be the length of the $i^{\scriptstyle th}$ row.  
We thus have $\nu_1\geq \nu_2\geq \cdots \nu_{Q_1}$. Once again, we 
restrict ourselves to the vacuum module, which is characterized by 
$\lambda_1-\lambda_2=0$ for even $L$. 

Let  us  call a {\em minimal} rigged configuration ($\nu,L,\{J=0\}$) 
with  all  the  integers  set  to  zero. This is the state with the 
smallest  value  of  the  momentum ascribed to the {\em excitations}. 
Let us evaluate what the sequence of ranks is under the ramification 
rules  that  map  this  rigged configuration into a standard tableau.  
Order the rows according to (\ref{inc}).

Note that $P_{\nu_1}=0$ for $L$ even and $P_j>0$ for $j<\nu_1$. 
Therefore $w_L=1$, and we delete the rightmost cell from the row           
$\nu_{m_{\nu_1}}$. The modified rigged configuration now has
a row of length $\nu_1 -1$ which is special by construction, while
$P_{\nu_1}\mapsto P_{\nu_1}+1$.  Therefore, $w_j=1,j=L,L-1,\ldots,
L-\nu_1+1$, until the top row has been {\em completely removed},
by which stage, $P_{\nu_2}=\nu_1$.  Before $P_{\nu_2}$ 
is brought back to $0$, $w_j=0,j=L-\nu_1,L-\nu_1-1,\ldots,L-2\nu_1+1$.
Thus, from a minimal configuration, strings of length $k$ generate
$k\, 0$'s and $k \, 1$'s in succession, and the longer strings 
are deleted first. 

In this way we end up with a word $w$:
\[
w=w_L w_{L-1} \cdots w_1 = 11\cdots 100\cdots 0
\underbrace{11\cdots 1}_{\nu_i}
\underbrace{00\cdots 0}_{\nu_i} 1\cdots 101\cdots10,
\]
where the first sequence of $1$'s and the first sequence of $0$'s
are each of length $\nu_1$, the next sequence of $0$'s and $1$'s of 
length $\nu_2$ each, and so on. In general, for a minimal configuration
as defined, there are sequences of 
$\nu_i \, 0$'s followed by $\nu_i \, 1$'s for $\nu\vdash {L\over 2}$.

\subsubsection*{Remark}

Given an admissible rigged configuration $(\nu,L,\{J\})$,  we define 
another, $(\nu,L,\{J^\sigma\})$ by performing the following operations 
(involutions) $\sigma_{j,\alpha}$ 
($j\in\{1,\ldots,\ell\}$,
$\alpha\in\{1,\ldots,m_j\}$) and 
$\sigma:=\prod_{j=1}^\ell \prod_{\alpha=1}^{m_j} \sigma_{j,\alpha}$: 
\[
\label{J-inv}
\begin{array}{rcl}
\sigma_{j,\mu}:\{J\}&\longrightarrow &\{J\}\\
J_{j,\mu}&\mapsto &J^\sigma_{j,\mu}=P_{j}-J_{j,m_j+1-\mu}.
\end{array}
\]
Note that $\sigma_{j,\alpha}^2=\sigma^2=$id.

Recall that a node
is added to a row in $\nu$ if it is the highest special row,
and this results in a decrease in the maximum integer allowed
for all rows above it. Every uninterupted sequence of say $x$
$1$'s would require the addition of cells to the {\em same} row, 
since by construction it would be assigned an integer equal 
to its maximum allowed, provided the integers rigging the rows
above it are less than the maximum allowed integer by at least $x$.
Thus every sequence of $(\nu_i)$ $1$'s would make a $\nu_i$-string 
pop up, but in order for the word of $0$'s and $1$'s to be a 
lattice permutation, every such sequence must be preceded by an
equal number of $0$'s.  Note also, that every row in the configuration
generated must be special by construction.

It is clear from the above discussion that $\sigma_{i,m_i}(\nu,L,\{J=0\})$
produces a word which can be obtained from the minimal configuration in the 
following way.  The first occuring sequence of $\nu_i \, 1$'s and
$\nu_i \, 0$'s ({\em i.e.}, that which starts at $w_k$ for the largest $k$)
is translated so as to start from $w_L$.  

\subsection{The maximal configurations:}

What about the maximal configuration, the one in which all the entries 
$\{J_{j,\alpha}=P_j\} \, \forall j$?  The previous remark makes it clear
that if $\sigma$ which is the product of {\em all} the involutions 
$\sigma_{j,\alpha}$ is applied, we would reach the maximal configuration,
which is a lattice permutation of the form
\[
w=w_L w_{L-1} \cdots w_1 = 
1010101\cdots 01100\cdots
\underbrace{11\cdots 1}_{\nu_i}
\underbrace{00\cdots 0}_{\nu_i} 
11\cdots 100\cdots 0,
\]
Thus, we conclude that the lattice permutations corresponding to 
$\{J=0\}$ and $\{J=P\}$,
$w_{\scriptstyle min}=w_1 w_2\ldots w_L$
and $w_{\scriptstyle max}=w'_1 w'_2 \ldots w'_L$ respectively, 
are related by the evacuation involution 
$w'_i=1-w_{L+1-i}$.

\subsection{The rest of the configurations}

Now we have to deal with the rest of the configurations $-$ we 
wish to show that any path, that is neither minimal nor
maximal, maps to a rigged configuration, that is also neither
minimal nor maximal, and {\em vice versa}.

More importantly, we want to show that if the ABF path has maximal height 
$\ell+ 1$, then the maximal string length is necessarily 
$\ell + 1$.  Let the length of the longest string be $l_{max}$.
If  at some stage, the height 
reached is $\ell+1$, then the corresponding partial tableau is such
that $\lambda_1 - \lambda_2 = \ell + 1$, and therefore
$P_{l_{max}}=\ell+1$ from equation 5.  

We know that for the rectangular tableaux (equal number of $0$ and $1$
occurences) which defines the vacuum sector, 
 $P_{l_{max}}=\lambda_1-\lambda_2=0$.  
Also, for every deletion from any row  
in $(\nu,\{J\})$, the corresponding $P_{l_{max}}$ goes 
{\em up} by $1$, and if the length $L$ is reduced by one, leaving
the shape of $\nu$ unchanged, all the 
$P_j$s would go {\em down} by $1$.
Thus, to go from $0$ to $\ell+1$, in order to 
reach the configuration described above, the number of 
deletions from the rows of $\nu$ must be $\ell+1$ more than 
reductions in $L$.  Every  time a row smaller than the longest 
is special, and cells 
have  to  deleted from it, the special string disappears quicker 
than the longest one would have (obviously). Therefore, the next 
step  just  after this disappearance would have to be a $0$, and the
length is reduced without deleting cells from $\nu$.  Therefore, the
only way $P_{l_{max}}$ can go from $0$ to 
deleted from a row of length at least $\ell+1$.  However, if the row has
more than $\ell+1$ cells, then the maximum height in the path picture
would then be more than $\ell+1$.

\subsection{Weighted sums}

In \cite{KR}, the following identity was established, based on the 
bijection described above:
\begin{equation}
\sum_{T\in SYT(\lambda)}q^{{\bf  p}(T)}=\sum_{T\in SYT(\lambda)} 
q^{{\bf c}(T)}=\sum_{\nu}\prod_{j=1}^{\ell} q^{c(\nu,L)} \biggl[ 
\begin{array}{c}
s_j+P_j\\
s_j
\end{array}\biggr].
\end{equation}
where
\begin{equation}
c(\nu,L)={1\over 2}L(L-1)-{1\over 2}LQ_1-\sum_{i}m_i P_i,
\end{equation}
and 
$$\biggl[ 
\begin{array}{c}
A+B\\
B
\end{array}\biggr]={\prod_{j=1}^{A+B} (1-q^j)\over
\prod_{j=1}^A (1-q^j)\prod_{j=1}^B (1-q^j)}.
$$
By the equality of the energy functional on paths, $H$ with 
${\bf \scriptstyle p}-{L\over 2}({L\over 2}-1)$, we arrive at a generating functional
for counting CTM paths as follows.                                  
Observe that $Q_j-Q_{j-1}=\sum_{k\geq j}m_k=-{1\over 2}(P_{j-1}-P_j)$,
so that 
$$m_j={L\over 2}\delta_{j1}-{1\over 2}\sum_{k}C_{jk}P_k$$
and thus, 
\begin{equation}
c(\nu,L)={1\over 2}L(L-1)-{1\over 4}L^2+{1\over 4}\sum_{ij}P_i C_{ij} P_j.
\end{equation}                                     
This establishes the equality of the
energy of the CTM path and the quasiparticle momentum of a rigged
configuration (\ref{qp}). 

In order to obtain the Virasoro character it is necessary to subtract off 
the Fermi momentum, or in the path language the energy of the ground state.
For the ground state, $P_j=0\, \forall j$ and $m_j=(L/2)\delta_{j1}$, so that 
$Q_1=(L/2)$.  After subtraction, we arrive at the proof of
\begin{equation}
B_L(\ell\Lambda_0)=\sum_{\nu}\prod_{j=1}^{\ell} 
q^{{1\over 4}\sum_{ij}P_i C_{ij} P_j} \biggl[ 
\begin{array}{c}
{1\over 2} L\delta_{j1}-{1\over 2}\sum_{k}C_{jk}P_k\\
P_j
\end{array}\biggr].
\end{equation}      

\subsubsection*{Remark}

The observation that Thomas' statistic ${\bf p}$ gives the required
energy function $H$ of a path ensures that the bijection of 
\cite{KKR,KR} will necessarily generalize the proof of the 
bosonic-fermionic identity for the higher rank cases.
(See \cite{forthcoming}.)  

\section{Discussion}

The point of this work was to establish an explicit weight-preserving 
bijection between two sets of combinatorial objects: The paths that 
parameterize the CTM eigenvectors, and the strings, or rigged 
configurations, that parametrize the RTM eihenvectors. In this bijection, 
standard tableaux play the crucial role of interpolating objects. The 
central point of this work, is that weights associated with each of the 
above objects: the energies of the paths, the momenta of the strings, 
and the Thomas statistic on the standard tableaux (which are related by
an involution to the Lascoux-Sch\"utzenberger charge), are equal.
Though our work is restricted to the vacuum sectors of the regime-III 
ABF models, this bijection is so natural, that we conjecture that 
analogous bijections exist for the rest of the sectors, and for
more general models.

\subsubsection*{Note added}

After this work was completed, we received a preprint \cite{NY},
where the relation between the energies of the CTM paths, and the
Lascoux-Sch\"utzenberger charge of tableaux was noted, and used
for somewhat different purposes.

\section*{Acknowlegements}

We wish to thank Paul Martin for discussions. S.~D. is supported 
by the EPSRC in the form of the grants GRJ25758 and GRJ29923. 
O.~F. is supported by the Australian Research Council. This work 
was partly done while S.~D. was visiting the University of Melbourne, 
and also while O.~F. was visiting the University of Utrecht. We 
thank the Australian Research Council, and the 
Netherlands Organization for Scientific Research (NWO) for finacial 
support of these visits.

\end{document}